\documentclass[12pt]{amsart}
\usepackage{color}
\usepackage{amssymb}

\usepackage{pdfpages}
\usepackage{graphicx}
\usepackage{dcolumn}

\RequirePackage[numbers]{natbib}
\RequirePackage{hyperref}
\usepackage{paralist}

\usepackage{amsmath}
\usepackage{fancyhdr}
\usepackage{amsthm}
\usepackage{amsfonts}

\usepackage{amscd}
\usepackage{latexsym}
\usepackage{graphics}
\usepackage{natbib}
\usepackage{afterpage}

\usepackage{subeqnarray} 
\usepackage{cases} 
\usepackage{booktabs}
 \usepackage{threeparttable}

\usepackage{subfigure}

\usepackage{float}

\graphicspath{{./figures/}}
\usepackage{amsfonts}
\usepackage{amsbsy}
\usepackage{fullpage}
\usepackage{natbib, mathrsfs}
\usepackage{verbatim}
\usepackage[latin1]{inputenc}
\usepackage{mhequ}
\usepackage{algorithm}
\numberwithin{equation}{section}

\newtheorem{theorem}{Theorem}[section]

\theoremstyle{plain}
\newtheorem{thm}{Theorem}
\newtheorem*{thm-non}{Theorem}

\newtheorem{example}[thm]{Example}

\theoremstyle{definition}
\newtheorem{defn}[theorem]{Definition}

\begin{document}

\title{Linear least square method for the computation of the mean first passage times of ergodic Markov chains}


\author{Yaming Chen}
\thanks{ Nanjing University Press, Nanjing, P.R.China. (cym\_shake@126.com)}


\maketitle







\begin{abstract}
~~An efficient and accurate iterative scheme for the computation of the mean first passage times ( MFPTs) of ergodic Markov chains has been presented. Firstly, the computation problem of MFPTs is transformed into a set of linear equations. It has been proven that each of these equations is compatible and their minimal norm solutions constitute MFPTs. A new presentation of the MFPTs is also derived. Using linear least square algorithms, some numerical examples compared with the finite algorithm of Hunter [LAA, 549(2018)100-122] and iterative algorithm of J.Xu [AMC, 250(2025)372-389] are given. These results show that the new algorithm is suitable for large sparse systems.
\end{abstract}



\section{Introduction}\label{I1}

 There are two kinds of algorithms for computing the mean first passage times(MFPTs for short) of finite ergodic Markov chains: finite algorithm and iterative algorithm.
 Meyer and Hunter propose finite algorithms for the computation of MFPTs based on generalized inverse matrix and stationary distribution vector.
In 1985, Grassman, Taksar and Heyman proposed the classical GTH algorithm and EGTH algorithm to find stationary distribution vectors. In 2016, Hunter gave rank-one update algorithm and perturbation method. Recent reviews of these finite algorithms are presented in Hunter's.

  In 2015, J. Xu first constructed an iterative method with parameters for the computing of MFPTs. Compared with the classical finite algorithm, the iterative method not only guarantees the accuracy of computation, but also greatly improves the scale and stability. At the same time, J. Xu also pointed out that the determination of iterative parameter $\alpha$ is a difficult problem. The parameter $\alpha$ is decided by the transition matrix and experimental calculation. In this paper, based on the definition equation of the mean first passage times matrix, the computation of MFPTs is transformed into a set of linear equations to be solved. These equations are all compatible, but the coefficient matrix may be singular. It is proved that the minimum norm solution of these equations form the mean first passage times matrix. Then, the least norm solution of these equations is solved uniformly by the linear least square algorithm (LS for short).
As far as the author knows, this is the first time for computing the MFPTs by LS algorithm. Numerical experiments show that this scheme is effective.

This article is organized as follows. In Section 2 the definition equation of the matrix of the mean first passage times of a $n$ state Markov chain and the related knowledge of the generalized inverse of the matrix are introduction without proofs. In Section 3  a set of linear equations are constructed, and the equivalence of this set of equations to the original matrix equation satisfied by MFPTs is proved. last section, A number of numerical examples are given. The LS algorithm is compared with Hunter's finite algorithm and Xu's iteration method in which the results show that the linear least squares algorithm is effective in computing the matrix of mean first passage times.

\section{Preliminaries}
  We set the scene by reintroducing the notation and theory that are statemented in \cite{H2017}. Let $\{X_n, n \geq 0\}$ be a finite  Markov chain with state space $S=\{1,2, \cdots, m\}$ and transition matrix $P=(p_{ij})$, where $p_{ij}=P\{X_n=j|X_{n-1}=i\}$ for all $i,j \in S $. In this article, we  focus on regular Markov chain(MC for short).
 the first passage time $T_{ij}$  is the length of time to go from a state $i$  to a state $j$  for the first time.
and define $m_{ij}=E[T_{ij}|X_0=i]$ as the mean first passage times from state $i$ to state $j$.
 The mean first passage matrix, denoted by $M=(m_{ij})$ , is the matrix with entries $m_{ij}=E[T_{ij}|X_0=i]$.

 \begin{thm} Let $P$ be a transition matrix of a homogeneous $n$ state Markov chain. $X_d$ be the diagonal matrix formed from the diagonal elements of matrix $X$. $J=[1]$(i.e. each element of $J$ is 1). The  matrix $M$ of mean first passage times is unique solution of the matrix equation\cite{K60,X15}
\begin{equation}\label{2a}
(I-P)X=J-PX_d.
\end{equation}
 \end{thm}

\begin{defn} \cite{Ben2003}  Let $A \in \mathbb{C}^{ m \times n}$. The Moore-Penrose generalized inverse of $A$ denoted by $A^+$  is the unique solution of four matrix equations:  $AXA=A, XAX=X, (AX)'=AX, (XA)'=XA$, where prime $'$ stands for the conjugate transposing of a matrix.
\end{defn}
   If $A \in \mathbb{C}^{ m \times m}$ and $A$ is nonsingular, then $A^+=A^{-1}$.

\section{Construction of a set of equivalent linear equations}
J. Xu\cite{X15} gives an equations with parameters $\alpha$ for MFPTs
\begin{equation}\label{2b}
(I-\alpha P)X_{k+1}=J+(1-\alpha)PX_k-PX_{kd},
\end{equation}
 which $\alpha$ satisfies $0\leq \alpha <1$. Then, Xu constructed an iterative algorithm with parameters $\alpha$ for computing MFPTs and proved the convergence of this algorithm\cite[Thm 3.2]{X15}.\\
\begin{algorithm}\mbox{(AlgXu)}
 Suppose that$(I-\alpha P)$admits an $LU$ factorization.

Choose an initial approximation $X^{(0)}\in\mathbb{R}^{n\times n}$, tolerance $\varepsilon>0$, and a dixed norm $\|\cdot\|$, starting with $k=1$.\\
(1)~Compute $Y^{(k-1)}=J+(1-\alpha)PX^{(k-1)}-PX_d^{(k-1)}$;\\
(2)~Solve $LY^{(k-1/2)}=Y^{(k-1)}$;\\
(3)~Solve $UX^{(k)}=Y^{\frac{k-1}{2}}$;\\
(4)~If $\|X^{(k)}-X^{(k-1)}\|<\varepsilon$, then return $M\approx X^{(k)}$ and stop; otherwise increase $k$ by 1 and continue with step (1).
\end{algorithm}
Compared with the classical finite algorithm\cite{H2017} , the algorithm AlgXu is of higher accuracy and better stability.
But it is still a difficult problem that choice of the best iteration parameter $\alpha$. In addition, the algorithm AlgXu has to solve two trigangular equations in every iteration step by step, and the iterative format (\ref{2b}) is not a standard iteration form, so it is difficult to analyze convergence.

Now we construct a set of equivalent linear equations to (\ref{2a}).

Let $e$ denote the column vector of $n$ ones and $X=(x_1, x_2, \cdots, x_n)$. Then the matrix equation (\ref{2a}) is equivalent to the following $n$ linear equations
\begin{equation}\label{2b1}
  x_i=Px_{id}+e,~~~~~~~~i=1, 2, \cdots, n,
\end{equation}
where $x_{id}=(x_{i1}, \cdots, x_{i,i-1}, 0, x_{i,i+1}, \cdots, x_{in})^T$.  Let $P_i$ be the matrix formed by replacing each  element of column $i$ of matrix P with 0 and  $A_i=(I-P_i)$. Then (\ref{2b1}) is equivalent to
\begin{equation}\label{2c}
  A_ix_i=e,~~~~~~~~i=1, 2, \cdots, n.
\end{equation}

The matrix $A=I-P$ is singular, but $A_i=I-P_i$ may be singular or nonsingular. However, it is always compatible for every $i$ to  equations (\ref{2c}). Therefore, whether $A_i$ is singular or not, its minimal norm solution is always $A_i^{+}e$.

The relation between the unique solution of (\ref{2a}) and the minimal norm solution of series equations (\ref{2c}) is given below.

\begin{thm}\label{2e}
Suppose the unique solution of matrix equation (\ref{2a}) is $X=(x_1,x_2, \cdots, x_n)$, then $x_i=A_i^{+}e, i=1, 2, \cdots, n$.
\end{thm}
Proof:  First, notice that for any matrix $Q$,  $QQ^+$  is a projection on $\mathbb{R}(Q)$ along  $\mathbb{N}(Q')$. If $x \in \mathbb{R}(Q)$,  then $QQ^+x=x$ or $(I-QQ^+)x=0$. \\
 \indent Next, let $P'_i=P-P_i$, $P=(p_1,p_2,\cdots,p_n)$,  $x_i=(x_{i1},x_{i2}, \cdots, x_{in})^T$, the minimal norm solution of of equations (\ref{2c}) is $x_i=A_i^{+}e$. Thus
$$(I-P)x_i=(A_i-P'_i)x_i=A_ix_i-P'_ix_i=A_iA_i^+e-P'_ix_i=e-x_{ii}p_i.$$
Therefore, if $X=(x_1,x_2,\cdots,x_n)$, then $(I-P)X=AX=J-PX_d$. $\Box$

In summary, the computation problem of MFPTs can be reduced to solving the minimal norm solution of a set of equations. As is known to all, there are many effective algorithms for solving the minimal norm solutions of systems with full rank or deficient rank.  Now, we use the linear least squares algorithm (LS for short)  based on QR decomposition  to compute $A_i^{+}e$. The implementation details of the LS algorithm  see \cite[Chapter 5]{GV96}, which is no longer given here.

\section{Numerical experiment and analysis}

 Two kinds of numerical experiments are done in this section. Finite algorithm and iterative method are chosen to compare with LS algorithm. Two finite algorithm are selected from the Hunter\cite{H2017}, the algorithm Proc2 and Proc4(HP2 and HP4 for short respectively), which have better computational efficiency for different scale problems, and the computing of the inverse of matrix $A$ in HP2 and HP4 is directly called the inv function command in program language. The iterative algorithm is selected from the J.Xu\cite{X15} algorithm.

 The three test index is selected from Hunter\cite{H2017}:

 1. Computational time. A comparison of results running on the same machine is given.

 2. Percentage of Zero errors(PZE). Let $M=\{m_{ij}\}$ be the mean first passage times, $\varepsilon_{ij}\equiv m_{ij}-\sum\limits_{k\neq j}p_{ik}m_{kj}-1$. PZE be the percentage of error terms $\varepsilon_{ij}$ that are zero.

3. Overall Residual errors(ORE). ORE$\equiv \sum\limits_{i=1}^{m}\sum\limits_{j=1}^m |\varepsilon_{ij}|$.

First, we compare three indicators above with 4 examples from different literatures, \cite{C15, H2017, X15}. In this part, each algorithm runs 20 times and then takes the average value. Secondly, we compare the advantages and disadvantages of finite algorithm and iterative algorithm with two large scale examples.

\begin{example}
In this case, we choose 4 test matrices $P_1$to $P_4$\cite{C15, H2017, X15} for comparison.
\begin{eqnarray*}
\fontsize{10}{12}\selectfont
 P_1=\left[
   \begin{array}{ccccc}
   0.136267 &0.292549 &0.266992 &0.220856 &0.083335\\
   0.198798 &0.019347 &0.129998 &0.321252 &0.330605\\
   0.246269 &0.215116 &0.044021 &0.249831 &0.244763\\
   0.400950 &0.149352 &0.012546 &0.303336 &0.133815\\
   0.200328 &0.084084 &0.351278 &0.337325 &0.026985
         \end{array}
         \right].
\end{eqnarray*}

\begin{eqnarray*}
\fontsize{10}{12}\selectfont
 P_2=\left[
   \begin{array}{cccccc}
0.268031 &0.255740 &0.201497 &0.265012 &0.007385 &0.002335\\
0.166582 &0.137728 &0.032748 &0.118446 &0.187835 &0.356660\\
0.093279 &0.226108 &0.081331 &0.206803 &0.094199 &0.298281\\
0.103853 &0.230590 &0.261709 &0.069110 &0.061473 &0.273265\\
0.101657 &0.261742 &0.128131 &0.002138 &0.204864 &0.301467\\
0.216100 &0.210158 &0.154059 &0.178624 &0.213131 &0.027928
         \end{array}
         \right].
\end{eqnarray*}

\begin{eqnarray*}
\fontsize{10}{12}\selectfont
 P_3=\left[
     \begin{array}{ccccc}
0.000000 & 0.701299& 0.298701 &0.000000 &0.000000\\
0.000000 & 0.000000 & 0.437907 & 0.562093 &0.000000\\
0.000000 &0.000000 & 0.000000 &0.632082 &0.367918\\
0.471475 &0.000000 & 0.000000 &0.000000 &0.528525\\
0.461323 &0.538677 & 0.000000 &0.000000 &0.000000
    \end{array}
    \right].
\end{eqnarray*}

\begin{eqnarray*}
\fontsize{10}{13}\selectfont
 P_4=\left[
     \begin{array}{ccccc}
0.999999 & 1e-7  & 2e-7  &  3e-7   & 4e-7 \\
0.4 &   0.3 &  0 &  0 &  0.3\\
 5e-7  & 0  & 0.999999  & 0  & 5e-7\\
     5e-7  & 0  & 0 &  0.999999  & 5e-7\\
     2e-7  & 3e-7  & 1e-7 & 4e-7  & 0.999999
    \end{array}
    \right].
\end{eqnarray*}
\end{example}

The compare results of the 4 examples above are shown in table \ref{L1},\ref{L2},\ref{L3}. In the first two examples, the PZE index of the LS algorithm is the best, followed by the two finite algorithm. In the last two examples, the PZE index of HP2 algorithm is the best, followed by LS and HP4A. Xu algorithm is not ideal.

\renewcommand{\arraystretch}{1} %
\begin{table}[H]
	\centering
	\setlength{\tabcolsep}{2.6mm} {
		\fontsize{10}{14}\selectfont
		\begin{threeparttable}
			\caption{Average of computation times (seconds)}
			\label{L1}
			\begin{tabular}{ccccc}
				\toprule
				& $P_1 $ & $P_2$ & $P_3$ & $P_4$   \cr
				\midrule
		HP2       &1.0312e-04 & .6912e-04 & 3.3034e-04   & 3.5310e-04  \cr
        XU        &0.0045 &  0.0068 &  0.0034    & 0.1450 \cr
		LS        & 0.0051 &  0.0060  &  0.0046   &    0.0069  \cr
						\bottomrule
			\end{tabular}
	\end{threeparttable} }
\end{table}

\renewcommand{\arraystretch}{1} %
\begin{table}[H]
	\centering
	\setlength{\tabcolsep}{2.6mm} {
		\fontsize{10}{14}\selectfont
		\begin{threeparttable}
			\caption{Percentage of Zero Error}
			\label{L2}
			\begin{tabular}{ccccc}
				\toprule
				& $P_1 $ & $P_2$ & $P_3$ & $P_4$   \cr
				\midrule
		HP2       &0.2000 & 0.2222 & 0.6400     & 0.6000  \cr
       	XU          &0.2900 & 0.1667 & 0    & --  \cr
		LS     & 0.4200 & 0.4722  & 0.5600     &  0.5300   \cr
						\bottomrule
			\end{tabular}
	\end{threeparttable} }
\end{table}

\renewcommand{\arraystretch}{1} %
\begin{table}[H]
	\centering
	\setlength{\tabcolsep}{3.6mm} {
		\fontsize{10}{14}\selectfont
		\begin{threeparttable}
			\caption{Overall Residual Error}
			\label{L3}
			\begin{tabular}{ccccc}
				\toprule
				& $P_1 $ & $P_2$ & $P_3$ & $P_4$  \cr
				\midrule
		HP2       &4.9652e-05 & 1.0995e-04 & 7.1054e-15  &7.5147e-09  \cr
        XU        &1.1976e-09 & 2.1252e-09 & 1.3786e-09   & --  \cr
		LS        & 1.3323e-14 & 2.2204e-14  & 7.1054e-15    &  7.6613e-09   \cr
						\bottomrule
			\end{tabular}
	\end{threeparttable} }
\end{table}

\begin{example}
\end{example}

Test matrix in this example is sparse irreducible transition matrix generated randomly \cite{H2017}.
The scale of the test matrix is from the 10 to the 510. The Xu algorithm does not work well, so it is no longer considered here. Only HP2, HP4A algorithm and linear least square algorithm with QR decomposition (LSQR) are compared.

The generation code for the test matrix $P$ is as follows\cite{H2017}£º\\
\indent Input n;\\
\indent  a=0.4;\\
\indent P=rand(n);\\
\indent P(P$>$a)=0;\\
 \indent  P=P-diag(diag(P))\\
\indent c=1./sum(P');\\
\indent for i=1:n \\
\indent \indent     P(i,:)=P(i,:)*c(i);\\
\indent end\\

In addition, the generated matrix  has been  checked for irreducibility.

\begin{figure}[H]
	\setlength{\abovecaptionskip}{0.cm}
	\setlength{\belowcaptionskip}{-0cm}
	\centering
         \includegraphics[width=12cm]{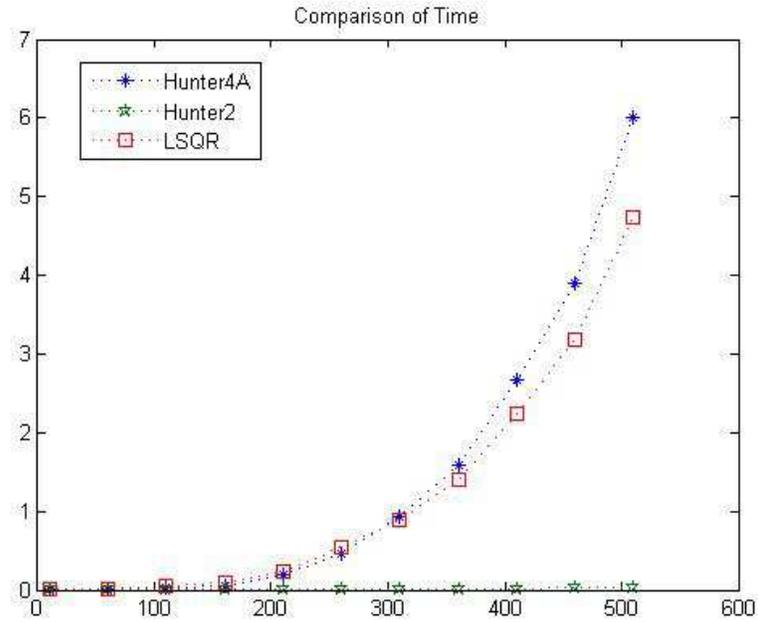}
	\centering
	\caption{Computation times for LS, HP2 and HP4A (random generating matrix).}\label{f3.1}
\end{figure}
\begin{figure}[H]
	\setlength{\abovecaptionskip}{0.cm}
	\setlength{\belowcaptionskip}{-0cm}
	\centering
		\includegraphics[width=11cm]{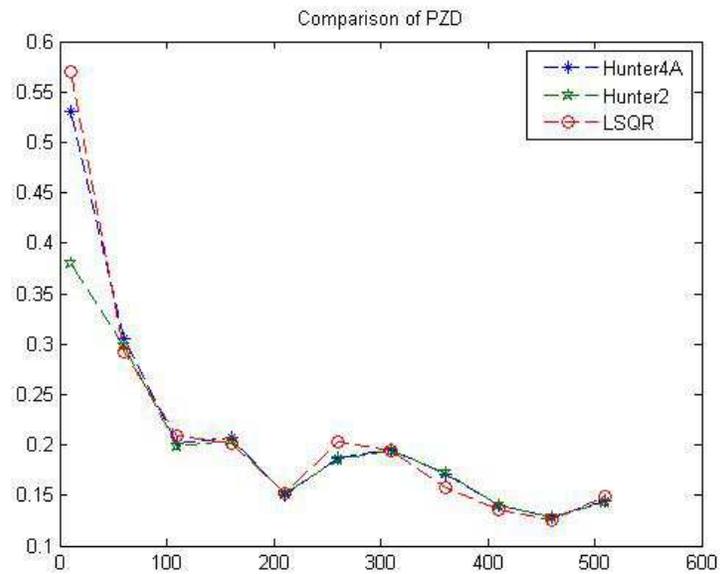}
	\centering
	\caption{PZE for LSQR, HP2 and HP4A (random generating matrix).}\label{f3.2}
\end{figure}

\begin{example}
\end{example}
Test matrix in this example is one dimensional random walk transition matrix\cite{C15}. The scale of the test matrix is from the 100 to the 2000. The efficiencies of HP4A algorithm and Xu algorithm are not good, so only HP2 and LSQR are compared.

Test matrix\cite{C15}
\begin{eqnarray*}
\fontsize{11}{14}\selectfont
P_8=\left[
     \begin{array}{cccccc}
0.75   &0.25   & 0      &\cdots &\cdots&0\\
0.25   &0.50   & 0.25   &\ddots &\ddots&\vdots\\
0      &0.25   & 0.50   &0.25   &\ddots&\vdots\\
0      &0      & \ddots &\ddots &\ddots&0\\
\vdots &\ddots & \ddots &0.25   &0.50  &0.25\\
0      &\cdots & \cdots &0      &0.25  &0.75
    \end{array}
    \right].
\end{eqnarray*}

\begin{figure}[H]
	\setlength{\abovecaptionskip}{0.cm}
	\setlength{\belowcaptionskip}{-0cm}
	\centering
         \includegraphics[width=12cm]{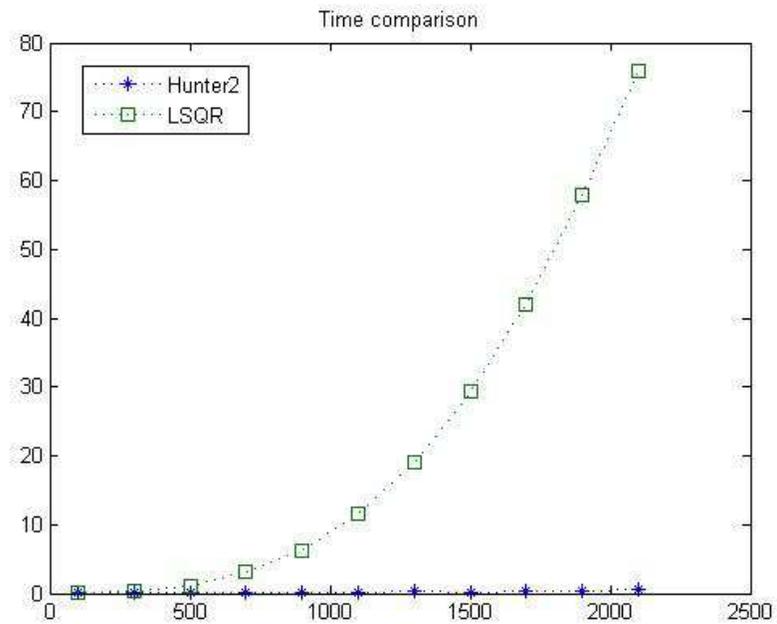}
	\centering
	\caption{Computation times for LSQR and HP2 (random walk matrix).}\label{f3.3}
\end{figure}

\begin{figure}[H]
	\setlength{\abovecaptionskip}{0.cm}
	\setlength{\belowcaptionskip}{-0cm}
	\centering
         \includegraphics[width=12cm]{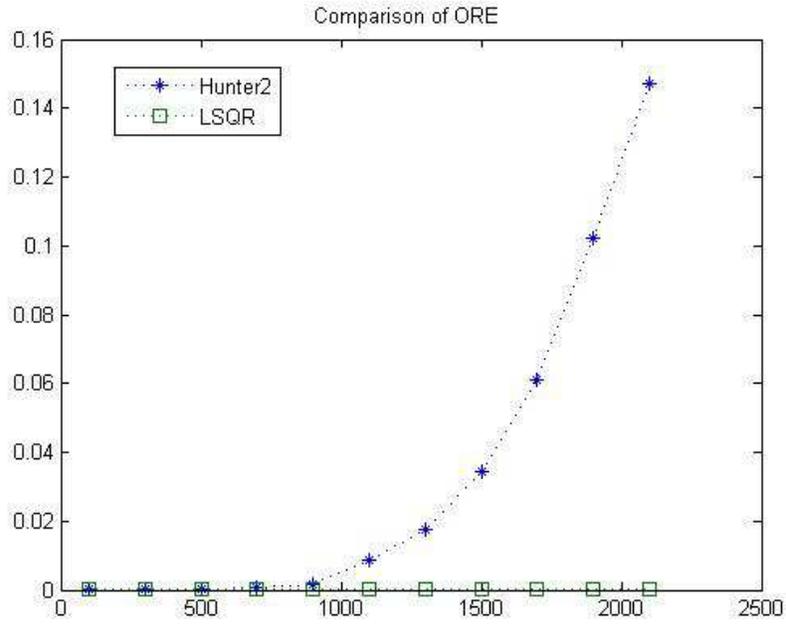}
	\centering
	\caption{ORE for LSQR and HP2 (random walk matrix).}\label{f3.4}
\end{figure}

\begin{figure}[H]
	\setlength{\abovecaptionskip}{0.cm}
	\setlength{\belowcaptionskip}{-0cm}
	\centering
         \includegraphics[width=12cm]{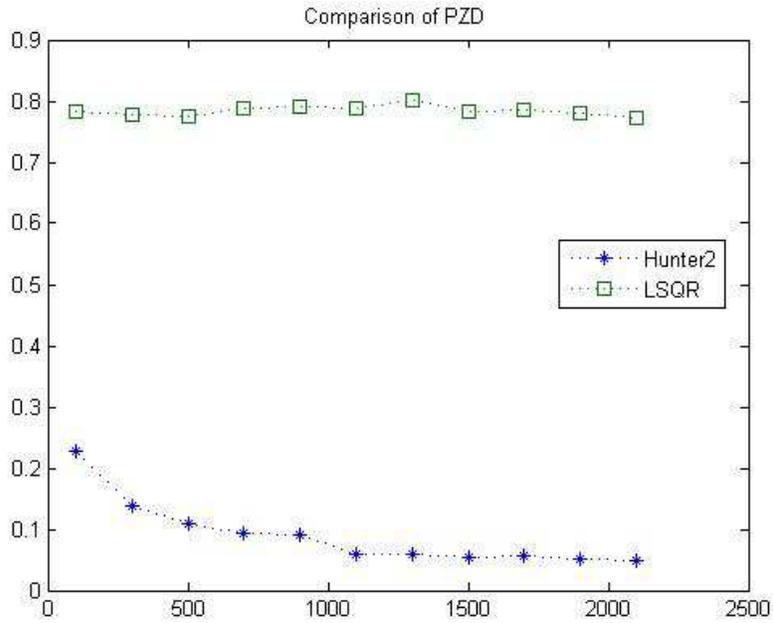}
	\centering
	\caption{PZE for LSQR and HP2 (random walk matrix).}\label{f3.5}
\end{figure}

From the results presented in Figs \ref{f3.3}, \ref{f3.4} and \ref{f3.5}, it is obvious that with the increase of matrix order, the computation time of LS algorithm increases greatly, while that of HP2 algorithm is relatively small. The index ORE of the LS algorithm is the magnitude of $1e-4$ and the HP2 algorithm magnitude is $1e-1$. For the important index, PZE, LS is stable around $0.8$, and HP2 decreases from $0.23$to $0.05$. It can be seen that LS algorithm can maintain high precision convergence in addition to long computation time.

\section{Conclusions}

The above examples show the effectiveness of the LS algorithm. But the computation time of LS algorithm is relatively long with finite algorithm HP2. It is noteworthy that $m$ independent linear equations are solved every time, so this kind of algorithm is very easy to implement in parallel. This may effectively reduce the computation time.

  In addition, people may be able to directly computation the generalized inverse matrix $A_i$ by parallel algorithm, so as to improve the computational efficiency.

\end{document}